\documentclass[12pt]{article}
\usepackage{amsfonts}
\begin{document}
\def\R{\mathbb{R}}
\def\C{\mathbb{C}}
\def\Z{\mathbb{Z}}
\def\N{\mathbb{N}}
\def\Q{\mathbb{Q}}
\def\D{\mathbb{D}}
\def\Ex{\mathbb{E}}
\def\P{\mathbb{P}}
\def\T{\mathbb{T}}
\def\Sp{\mathbb{S}}
\def\hb{\hfil \break}
\def\ni{\noindent}
\def\i{\indent}
\def\a{\alpha}
\def\b{\beta}
\def\e{\epsilon}
\def\d{\delta}
\def\De{\Delta}
\def\g{\gamma}
\def\qq{\qquad}
\def\L{\Lambda}
\def\E{\cal E}
\def\G{\Gamma}
\def\F{\cal F}
\def\Hi{\cal H}
\def\K{\cal K}
\def\O{\cal O}
\def\A{\cal A}
\def\B{\cal B}
\def\L{\cal L}
\def\M{\cal M}
\def\N{\cal N}
\def\Om{\Omega}
\def\om{\omega}
\def\s{\sigma}
\def\t{\theta}
\def\th{\theta}
\def\Th{\Theta}
\def\z{\zeta}
\def\p{\phi}
\def\m{\mu}
\def\n{\nu}
\def\l{\lambda}
\def\Si{\Sigma}
\def\q{\quad}
\def\qq{\qquad}
\def\half{\frac{1}{2}}
\def\hb{\hfil \break}
\def\half{\frac{1}{2}}
\def\pa{\partial}
\def\hb{\hfil \break}
\def\ni{\noindent}
\def\i{\indent}
\def\half{{1 \over 2}}
\def\Om{\Omega}
\def\om{\omega}

\begin{center}
{\bf PREDICTION THEORY FOR STATIONARY FUNCTIONAL TIME SERIES}
\end{center}
\begin{center}
{\bf N. H. BINGHAM}
\end{center}

\ni {\it Abstract} \\
\i We survey aspects of prediction theory in infinitely many dimensions, with a view to the theory and applications of functional time series. \\

\ni {\it Keywords} \ Cram\'er representation, Kolmogorov isomorphism theorem, Verblunsky coefficients, Szeg\H{o}'s theorem, Szeg\H{o} alternative, Beurling-Lax-Halmos theorem, functional time series, functional principal components, Karhunen-Lo\`eve expansion, kernel methods \\

\ni 2010 Mathematics Subject Classification: Primary 60-02 Secondary 62-02\\

\ni {\bf 1. Introduction} \\ 

\i This paper continues the theme of the author's earlier surveys [Bin1] on prediction theory for one-dimensional time series, [Bin2] on the finite-dimensional case, and (with Badr Missaoui, [BinM]) touching briefly on the infinite-dimensional case, our theme here.  Our motivation is partly mathematical interest and completeness, partly the vigorous development of functional data analysis (FDA; [RamS1], [RamS2], [HorK]) made possible by the explosive growth in computer power, data storage and data handling.  \\
\i We begin in \S 2 with the Cram\'er Representation (CR) and the Kolmogorov Isomorphism Theorem (KIT), on which everything rests.  In \S 3 we turn to Verblunsky coefficients, Szeg\H{o}'s theorem and the Wold decomposition, followed in \S 4 by the Szeg\H{o} alternative and factorization, all themes familiar from e.g. [Bin1], [Bin2].  Section 5 is on the Beurling-Lax-Halmos theorem and inner functions.  Section 6 is on numerical implementation.  Complements follow in \S 7, in particular (\S 7.1) the deterministic case and (\S 7.2) model spaces.  We close in \S 8 with some open questions, with which the area abounds. \\
\i The mathematics here involves vectorial integration (see e.g. [Rud2, Ch. 12]).  As with FDA, implementation involves discretization (`calculus is continuous, calculation is discrete'), and so in principle (subject to a suitable choice of dimension) reduces the setting to a finite-dimensional one, which can be handled by the finite-dimensional methods of e.g. [Bin2] and the references there.  We prefer to use the infinite-dimensional context suggested by the nature of the data.  See e.g. \S 6.1 below, [RamS1], [DetKA] for more on this.  (Ramsay and Silverman remark [RamS1, p.11] `In general, prediction theory is beyond our scope, and is only considered here and there'.) \\        

\ni {\bf 2.  The Cram\'er Representation and the Kolmogorov Isomorphism Theorem} \\

\ni {\it 1.  The Cram\'er Representation} \\
\i We confine ourselves throughout to {\it stationary} processes, that is, those whose distributions are invariant under {\it time-shifts}.  One can work in discrete or continuous time, depending on preference, context or the data available.  As the term `time series' indicates, the first is the traditional one, and we shall follow it here.  Then the time set is the integers $\Z$, with (Pontryagin-van Kampen) dual [Pon, Ch. 6] the torus $\T$ (equivalently $\R/2 \pi \Z$; one can pass between these by $\th \leftrightarrow e^{i \th}$, and we shall do this at will).  In the second case, the time-set is the real line $\R$, with dual $\R$ also.  Stationarity is a strong condition, which yields correspondingly strong conclusions; we turn later to how it may be relaxed (\S 7.4).  For a monograph treatment, see Nikolskii [Nik1], whose sub-title `Spectral function theory' gives a hint of the mathematics involved (see also [Nik2]).    \\
\i We write our process (or time series) as $x = \{ x_n: n \in \Z \}$, where the $x_n$ are in $\R$, ${\R}^d$, a Hilbert space $\Hi$ or a Banach space $\B$, depending on context (one can work more generally; see e.g. [BinM, \S 5.5, \S 5.6]).  While our data consists of functions, as in our title, these functions will always belong to function spaces, which are at least topological vector spaces, the elements of which we will call vectors as usual.  So vectors here are {\it infinite-dimensional} unless otherwise stated.  Covariances are matrices in the finite-dimensional case (as in multivariate analysis in statistics, and in [Bin2]).  But here, covariances are (linear) {\it operators}, hence the crucial role of operator theory in what follows.  We will work in a Hilbert space $\Hi$ unless otherwise stated. \\
\i Write the time-shift $n \mapsto n+1$ as $U$.  Then $U$ is unitary, and generates a unitary group, $\cal U$.  Being unitary, $U$ is normal, and so the spectral theorem for the unitary case (Stone's theorem) applies ([Rud2, Ch. 12]; [Sto VIII.2], [RieN, \S 109], [DunS, X.2]).  This involves integration over $\T$ with respect to a {\it resolution of the identity} (projection-valued measure, in Mackey's terminology [Mac]) $E = \{ E(\th): \th \in \T \}$: $E(\T) = 1$, and   \\
$$
U^n = \int_{\T} e^{i n \th} dE(\th), \qquad n \in \Z          \eqno(U)
$$
(as all our integrals will be over the torus $\T$, we omit the $\T$ below).  Note that as $|e^{in \th}| = 1$ and $E(\T) = 1$, for $x \in \Hi$ $(U)$ gives 
$$
\Vert U x \Vert = \Vert \int e^{i \th} dE(\th) x \Vert \leq \Vert x \Vert:
$$
$U$ is a {\it contraction}.\\
\i By above, 
$$
x_n = U^n x_0 = \int e^{i n \th} dE(\th) x_0
               = \int e^{i n \th} dY(\th), \qquad n \in \Z,  
$$
giving the {\it Cram\'er representation}
$$
x_n = \int e^{i n \th} dY(\th), \qquad n \in \Z.                    \eqno(CR)  
$$ 
The random measure $Y$ here is the {\it Cram\'er measure} [Cra1], [Cra2], [CraL], or (Cram\'er's terminology) the {\it spectral process}, of the stationary process $x$.  This is {\it orthogonally scattered} (Masani, [Mas1]): the masses of disjoint sets are orthogonal (see also Aue and van Delft [AuevD], Prop. 2.1). \\

\ni {\it 2. The Kolmogorov Isomorphism Theorem} \\
\i One has the {\it Kolmogorov Isomorphism Theorem (KIT)} ([Kak3, p.104]; [ManS]; [vDelE])]: as in the scalar and matrix cases (see e.g. [Bin1], [Bin2]),
$$
x(n) \leftrightarrow e^{in.} \ I, \qquad n \in \Z,         \eqno(KIT)
$$
with $I$ the identity operator (here the {\it time domain} is on the left, the {\it frequency domain} on the right). \\
\i We recall Chung's famous dictum ``{\it The process is the thing}".  A stochastic process is an infinite-dimensional object, characterised by its {\it distribution}.  Both $(KIT)$ and $(CR)$ deal with the distribution of the process in full.  They are really two different aspects of the same thing, with the distribution implicit in $(KIT)$ and explicit in $(CR)$.  By contrast, we shall meet below results that deal only with the {\it second-order aspects} of the process: (mean and) covariance.  The key result here is {\it Verblunsky's theorem}, or the {\it Verblunsky isomorphism} (\S 3), characterising (or parametrising) the covariance structure.  This is fully informative about the process when the process is {\it Gaussian}; we discuss Gaussianity in \S 2.4 below (cf. [Bin4, \S 4.1] in continuous time). \\

\ni {\it 3.  The Gramian} \\  
\i Take $\Hi$ a separable Hilbert space, $B(\Hi)$ the algebra of all bounded linear operators on $\Hi$, $T(\Hi)$ the trace class operators in $B(\Hi)$.  The Kolmogorov Isomorphism Theorem in the Hilbert case is due to Mandrekar and Salehi [ManS]; here we follow Kakihara [Kak1], [Kak2], [Kak3], [Kak4].  Write $X := L_0^2(\Om, \Hi)$ for the Hilbert space of all $\Hi$-valued (strong) random variables with mean 0 and finite second moment: with $(\Om, {\cal F}, \P)$ the probability space, $\Ex [x] := \int_{\Om} x(\om) d \P(\om) = 0$, 
$$ 
\Vert x \Vert^2_X := \int_{\Om} \Vert x(\om) \Vert^2_{\Hi}\  d \P (\om) 
= \Ex \ [ \ \Vert x(\om) \Vert^2_{\Hi}] < \infty,
$$  
and so with inner product on $X$
$$
(x,y)_X := \int_{\Om} (x(\om), y(\om))_{\Hi} \ d \P(\om) 
= \Ex [(x,y)_{\Hi}].
$$
The {\it Gramian operator}, or {\it Gramian}, is the $T(\Hi)$-valued inner product $[.,.]_X$ defined by
$$
([x,y]_X \phi, \psi)_H 
= \int_{\Om} (x(\om), \phi)_{\Hi} (\psi, y(\om))_{\Hi} \ d \P(\om)
$$
$$
{ \ \ \ \ \ \ \ \ \ \ }
= \Ex [(x,\phi)_{\Hi} (\psi,y)_{\Hi}] \quad (x,y \in X, \phi, \psi \in 
\Hi).                                                          \eqno(GO)
$$
In the finite-dimensional case, this reduces to a doubly-indexed set of scalars, which gives a matrix, the Gramian matrix [HorJ, \S 7.2], which is positive definite; so too is the Gramian operator likewise.  Symbolically,
$$
[x,y]_X = \int_{\Om} x(\om) \otimes \overline{y(\om)} d\P(\om) 
= \Ex[x \otimes \overline{y}],
$$
where
$$
(\phi \otimes \overline{\psi}) {\phi}' 
:= ({\phi}', \psi)_{\Hi} \phi, \qquad \phi, {\phi}', \psi \in {\Hi}.
$$
Now $X$ is both \\
(i) a left $B(\Hi)$-module under the module action
$$
(a,x) \mapsto a.x = ax \qquad (a \in B({\Hi}), x \in {\Hi}),
$$
(ii) a Hilbert space with Gramian $[.,.]_X$: \\
$X$ is a {\it normal Hilbert} $B(\Hi)$-{\it module} in the sense of [Kak2]. (For background on Hilbert modules, see Lance [Lan].  Both Hilbert modules and Gramians go back here at least to Masani's comments on Wiener's work on prediction [Mas2].) \\
\i We may now without ambiguity abbreviate $[.,.]_X$ (whose values are non-random operators) to $[.,.]$.  With it, we can define the {\it spectral measure}, or {\it control measure},
$$
F(.) := [Y(.), Y(.)]                                     \eqno(SM)
$$
(operator-valued, indeed Gramian-valued), and the {\it operator covariance function}, 
$$
\G(m,n) := [x(m), x(n)].
$$
The process $x = \{ x(n) \}$ is called {\it operator stationary} (or just {\it stationary}) if its operator covariance function (or just {\it covariance function}) is a function of $m-n$ only, $\tilde \G(m - n) = \G(m,n)$, say.  We assume stationarity unless otherwise stated. \\
\i As in the matrix case above, the Gramian is {\it positive definite} (Gretsky [Gret, \S 3]).  One has [Kak3] the {\it spectral representation}
$$
\tilde \G(n) = \int e^{i n \th} dF(\th).                  \eqno(SR)
$$
This is the operator version of Herglotz's theorem (or Bochner's theorem); cf. [Bin1], [Bin2] in finite dimensions and van Delft and Eichler [vDelE] in this setting. \\  
\i Below, we shall need the {\it spectral density} $f$, the Radon-Nikodym derivative of the absolutely continuous component of the spectral measure $F$ with respect to normalised Lebesgue measure $d \th / 2 \pi$ on $\T$, operator-valued, as $F$ is (take $f$ as $0$ if $F$ is singular).\\

\ni {\it 4. Gaussianity} \\
\i In $(CR)$, one has the stationary process $x = (x_n)$ in the time domain represented as a stochastic integral of the orthogonal-increments process $Y = (Y(\theta))$ in the frequency domain.  For Gaussians, orthogonality is the same as independence, so if $Y$ is Gaussian, it has independent increments.  As sums (and so integrals) of independent Gaussians are Gaussian, $x$ is then Gaussian by $(CR)$.  The converse also holds (see e.g. [Cra3], [Horo]).  All this is true {\it regardless of the dimension} $d$. \\
\i This simple and basic fact should thus have been stated (indeed, stressed) in the $d=1$ case in [Bin1] and the case $1 < d < \infty$ in [Bin2].  Unfortunately, Gaussianity is only touched on in [Bin1]  (\S 2, KIT, \S 4, Rajchman measures, \$ 6.1, $\phi$-mixing), and not even mentioned in [Bin2].\\ 
  
\ni {\bf 3.  Verblunsky coefficients; Szeg\H{o}'s theorem; Wold decomposition} \\

\ni {\it 1. Verblunsky's theorem} \\
\i In the scalar case [Bin1], the distribution of a (discrete-time, complex-valued) stationary sequence may be described (encoded, parametrised) by a sequence $\a = ({\a}_n)_0^{\infty}$ of complex numbers ${\a}_n \in \D$, the unit disc of $\C$, the {\it Verblunsky doefficients} (there are several other names; see [Bin1], [Sim1]). This is {\sl Verblunsky's theorem}.  The relevant theory depends heavily on orthogonal polynomials on the unit circle (OPUC), due originally to Szeg\H{o}, the theme of Simon's books [Sim1], [Sim2] ([Sim3] deals with orthogonal polynomials on the line and the cirle together; these correspond to continuous and to discrete time).  \\
\i In the finite-dimensional ($\ell$-vector, $\ell \times \ell$ matrix) case, the relevant theory is matrix orthogonal polynomials on the unit circle (MOPUC) [Bin2].  The process also has a Verblunsky parametrisation, but now the ${\a}_n$ are ($\ell \times \ell$) matrices, of norm $\Vert \a_n \Vert < 1$.  These encode the stationary processes, as before: {\it Verblunsky's theorem}.  See Damanik, Pushnitski and Simon [DamPS, (3.10), Th. 3.12]. \\

\ni {\it 2. Szeg\H{o}'s theorem} \\ 
\i For {\it stationary} processes, the key result in both cases is {\it Szeg\H{o}'s theorem}, which relates to the influence of the {\it remote past}.  This may be absent, e.g., bathwater forgetting its thermal history as it thermalises; total, e.g. tempered steel, whose thermal history is locked in; or present with a partial influence, e.g. the climate-weather interplay, where (on a time-scale, in years, short enough to neglect climate change) climate is permanent, while weather is (again on a suitable time-scale, in days) temporary and unpredictable, indeed chaotic. \\ 
\i Szeg\H{o}'s theorem in the scalar case tells us when the influence of the remote past is {\it not total}.  The condition for this is {\it Szeg\H{o}'s condition} 
$$
\log f \in L_1(\T).                                              \eqno(Sz)
$$  
This is extended to the finite-dimensional setting by Derevyagin, Holtz, Khrushchev and Tyaglov [DerHKT, Th. 28, Th. 29]: with $\dag$ for the adjoint (following their notation here), det and tr for determinant and trace,
$$
\prod_0^{\infty} \det (1 - {\a}_k {\a}_k^{\dag})
= \exp \int \hbox{tr} \ \log f(\th) d \th /2 \pi              \eqno(KSz)
$$
(the {\it Kolmogorov-Szeg\H{o} formula}; see e.g. [Bin1, \S 4]), and Szeg\H{o}'s condition -- that the right here is positive -- holds iff
$$
\sum_0^{\infty} \Vert {\a}_k^{\dag} {\a}_k \Vert < \infty.
$$
(extending earlier results of Delsarte, Genin and Kamp [DelGK, Th. 18, 19]).   Note that $(KSz)$ involves the {\it determinants} of the matrices ${\a}_n$, and (despite the Fuglede-Kadison determinant, [BleL]) these do not extend to the infinite-dimensional case in general.  \\
\i Payen [Pay] makes a thorough study of the Hilbert-valued case.  There ([Pay, II.]; [BinM]) he gives an infinite-dimensional form of Szeg\H{o}'s theorem in terms of factorization (\S 4 below). \\
\i Call the last display 
$$
\a \in {\ell}_2(\mathbb{N}).                           \eqno(\a - {\ell}_2)
$$  
The corresponding result for  
$$
\a \in {\ell}_1(\mathbb{N})                            \eqno(\a - {\ell}_1)
$$ 
is {\it Baxter's theorem} ([Sim1, Ch. 5]; [Bin1, \S 5] for the scalar case, extended by Kasahara and Bingham [KasB1, \S 5]).  A thorough study of Baxter's theorem in the matrix case was given by Dym and Kimsey [DymK]; cf. [Bin2], [KasB3]. \\  

\ni {\it 3. Wold decomposition} \\
\i The {\it past at time} $n$, the {\it remote past} and the {\it remote future} of the process $x = (x(n))$ are the closed linear subspaces spanned by the random variables below: 
$$
 H(x, n) := \bigvee (x(k): k < n), \quad 
 H(x, - \infty) := \bigcap \ (\bigvee ( x, n): n \in \Z ),
 $$
 $$
 H(x, + \infty) := \bigcap \ (\bigvee ( x, n): n \in \Z ).  
 $$
 The process is called {\it deterministic} if all three are equal, {\it purely non-deterministic (pnd)} if the second is trivial.  The {\it Wold decomposition} [Kak2] splits the process $x = (x(n))$ into a deterministic and a purely non-deterministic component,
 $$
 x = x_d + x_p,
 $$
 which are Gramian orthogonal:
 $$
 [x_d(m), x_p(n)] = 0 \quad (m, n \in \Z).
 $$
 \i The Cram\'er measure $Y$ of the process also splits, into a sum of absolutely continuous and singular components, $Y_a$ and $Y_s$, and similarly for the spectral measure:
 $$
 Y = Y_s + Y_a, \qquad F = F_s + F_a.
 $$
 In one dimension, one has {\it Wold-Cram\'er concordance}: regarded as processes in their own right, the `good' process $x_p$ has Cram\'er and spectral measures $Y_a$, $F_a$, while the `bad' one $x_d$ has $Y_s$, $F_s$ (see [Bin1] for details and references).  In the matrix case, Wold-Cram\'er concordance holds in the full-rank case but not in general ([Bin2]; see e.g. Payen [Pay, Remarque 8, 376-7]).  The vector case is studied in [Kak3], where such concordance is shown to be preserved under dilation from the stationary case here to the harmonisable case ($\S$ 7.4). \\
 \i The `good' component $x_p$ in the Wold decomposition is a {\it moving average} of the products of the {\it innovations} (new randomness) at past times and the matrices appearing in the Taylor expansion of the Szeg\H{o} function (`analytic square root' of the spectral density matrix), as we shall see below.  Compare [Bin4, \S 2.5] in continuous time. \\
 \i In the Hilbert-valued case, spectral criteria for the process to be purely non-deterministic are given by Kallianpur and Mandrekar [KalM]. \\    
  
\ni {\bf 4.  Szeg\H{o} alternative; factorization} \\

\ni {\it 1. Szeg\H{o} alternative} \\
\i In one dimension, one has a clean split, the {\it Szeg\H{o} alternative}, between `good' and `bad' cases.  In the first, there is a genuine innovation (input of new randomness) in each time-step from $n$ to $n+1$.  The size of this new input is measured by the {\it prediction error variance} ${\s}^2$, the infimum of the variances of all linear predictors based on the present and (finite sections of) the past of the next future value $x_{n+1}$.  This `good' case can only happen if there is an absolutely continuous component $F_a$ to $F$, that is, if the spectral density $f$ is not a.e. zero.  When this is so, it happens if and only if the {\it Szeg\H{o} condition} $(Sz)$ holds, and then ${\s}^2$ is the geometric mean of $\log f$ (see e.g. [Bin1]), so $(KSz)$ becomes
$$
{\s}^2 = \exp \{ \int \log f(\th) d\th / 2 \pi \} > 0.          
$$
Note that the singular component of $F$ (if present) plays no role here. \\
\i In the Wold decomposition, $x$ is the sum of a `bad' part $x_d$ ($\S 3$), and the `good' part $x_p$, a moving average of the innovations over the past to date; this contains a factor $\s$, and so is absent if $\s = 0$ [Bin1].  In terms of the Verblunsky coefficients,
$$
{\s}^2 = {\prod}_1^{\infty} (1 - |{\a}_n|^2),
$$
so $(Sz)$ holds iff $\s > 0$ iff the product converges, i.e. iff $(\a - {\ell}_2)$ holds.  Then (and only then), one can define the {\it Szeg\H{o} function}
$$
h(z) := 
\exp \Bigl( \frac{1}{2} \int \Bigl(\frac{e^{i\t} + z}{e^{i\t} - z} \Bigr) \log f(\t) \ d\t / 2 \pi \Bigr) \qquad (z \in D). \eqno(OF)
$$
This is an {\it outer function}, in the Hardy space $H^2(\D)$ ([Dur]; [Sim1]; [Bin1]).  It is analytic in $\D$, and zero-free there [Sim1, Th. 2.4.1]. As (from Fatou's theorem: [GarMR, Th. 1.10], [Rud1, Th. 17.10])
$$
 \int \Bigl(\frac{e^{i\t} + z}{e^{i\t} - z} \Bigr) \log w(\t) \ d\t / 2 \pi \to \log f(\phi) \ \ a.e. \ \ \quad (z = r e^{i \phi}, \quad r \uparrow 1),
 $$
its square has radial limit $f$ a.e. on $\T$ :
$$
\lim |h(z)|^2 = f(\th) \ \ a.e. \ \ (z = r e^{i \th}, \ r \uparrow 1).
$$
It may thus be regarded as the {\it analytic square root} of the spectral density $f$. \\
\i In the matrix case, we refer to [Bin2, \S 5] for details and references.     Szeg\H{o}'s theorem (the matrix forms of $(Sz)$, $(KSz)$ above) holds; see [DerHKT, Th. 28].  As noted above, this involves {\it determinants}, which do not extend to the infinite-dimensional case. \\

\ni {\it 2. Factorization} \\
\i In the scalar case, existence of the Szeg\H{o} function (analytic square root of the spectral density), the Szeg\H{o} condition $(Sz)$, and $\a \in {\ell}_2(\mathbb{N})$ (above) are equivalent. One seeks a matrix and a vector version of this.  So one seeks to factorise a matrix spectral density into its analytic square root times its adjoint (denoted by $\ast$ here, to conform to the sources below),
$$
F = \Phi {\Phi}^*.
$$
There is an extensive theory here, due mainly to Wiener and Masani and to Helson and Lowdenslager.  We note that the coefficient matrices ${\Phi}_n$ in the Taylor expansion
$$
\Phi (z) = \sum_0^{\infty} {\Phi}_n z^n
$$
appear in the prediction-error matrix (and operator, in the vector case, below).  See e.g. Masani [Mas2, p.278], Whittle [Whi] (for the scalar case, see e.g. [GreSz, \S 10.8]).   \\
\i For factorizations in the operator case, we refer to [RosR Ch. 6].  We also mention briefly the approach of Power [Pow].  We make the restrictive assumption that the spectral density $f$ is essentially bounded (i.e. is bounded, after excluding some null set on the torus $\T$).  Then
$$
f = h h^* + g,
$$
where $g$ is a positive operator and $h$ is analytic and outer.  Here $g$ can be taken to be {\it minimal}.  A Wold(-Zasuhin) decomposition is obtained.  The prediction-error operator $G(f)$ of $f$ is obtained as
$$
G(f) = G(h h^*) = (QhQ)(QhQ)*,
$$
where $Q$ is an orthogonal projection; we refer to [Pow] for details.\\    
\i As Power points out, his method gives $g$ and $h$ as functions of $f$, while the methods of Sz.-Nagy and Foia\c s ([SzNFBK]; \S 5 below) do not.  But this applies only to the (`good' and) `bad' behaviour encoded in the spectral {\it density}, the absolutely continuous component in the Lebesgue decomposition of the spectral {\it measure}.  We turn now to how to address the `bad' behaviour encoded in the other two components. \\
\\
\\      
          
\ni {\bf 5.  The Beurling-Lax-Halmos theorem;  inner functions} \\

\i So far as prediction theory is concerned, it is the `good' part of the process (or its spectral measure) that matters.  Nevertheless, the Wold decomposition shows us that to understand the structure of a stationary process we need to look at the `bad' part also. \\
\i Recall the {\it outer function} of $(OF)$, \S 4 above.  The term is due to Beurling [Beu], as is the corresponding term {\it inner function}.  These correspond respectively to the `good' and `bad' parts above.  When we pass from the time domain (Wold decomposition) to the frequency (or spectral) domain by the Kolmogorov Isomorphism Theorem, we obtain the {\it factorization} into outer and inner factors in Hardy spaces $H^2(\Hi)$.  For the scalar case (where the Hardy space is written $H^2$) there are good accounts in Duren [Dur], Garnett [Garn] and Hoffman [Hof] (see also [Rud1, Ch. 17]).  For the vector case, see [SzNFBK, III, V], [Nik1, I.7, XI.3], [RosR, Ch. 4-6], [Hof, 114-116]. \\
\i The remote past (the `bad' part above) is invariant under the time-shift.  The invariant subspaces are exactly those given by multiplication in $H^2(\Hi)$ by an inner function ([SzNFBK, V, Th. 3.3]; [Nik1, I.7, XI.3]; [RosR, 1.12]).  This is the Beurling-Lax-Halmos theorem, due to Beurling (dimension $d = 1$), Lax [Lax] ($1 < d < \infty$) and Halmos [Hal] ($d \leq \infty$).  So study of the remote past reduces to study of inner functions, $u$ say.  These have a rich structure; they factorise into a unimodular constant, a Blaschke product, and an integral factor as in $(OF)$ but with a singular rather than an absolutely continuous measure.  See \S 7.1,2 below.  \\
\i Halmos's approach [Hal] uses the concept of a {\it wandering subspace}; these correspond to the innovations (new randomness) in the process.  For an isometry $V$ on a Hilbert space $\Hi$, call a subspace $L$ of $\Hi$ {\it wandering} if for distinct integers $m, n$ $V^m L$ and $V^n L$ are orthogonal.  Then [SzNFBK, \S 1.1] if
$$
M_+(L) := {\oplus}_0^{\infty} \ V^n L,
$$
one has
$$
L = M_+(L) \ominus V M_+(L).
$$  
Call an isometry $V$ on $\Hi$ (which will be $U$ of \S 2 for us)  a {\it unilateral shift} if $\Hi$ has a wandering subspace $L$ with $M_+(L) = \Hi$.  Then $L$, called {\it generating} for $V$, is uniquely determined by $V$, indeed $L = \Hi \ominus V \Hi$.  One has the {\it Wold decomposition} ([SzNBK, Th. 1.1], [RosR, \S 1.3]): $\Hi$ decomposes into the orthogonal sum ${\Hi} = {{\Hi}_0} \oplus {{\Hi}_1}$ such that the ${\Hi}_i$ reduce ${\Hi}$ (each is mapped onto itself by $V$), $V | {\Hi}_0$ is unitary and $V | {\Hi}_1$ is a unilateral shift.  This decomposition is unique; indeed,
$$
{{\Hi}_0} = \bigcap_0^{\infty} V^n {\Hi}, \quad {{\Hi}_1} = M_+(L), \ \hbox{where} \ L = {\Hi} \ominus V {\Hi}.
$$
Either ${{\Hi}_0}$ or ${{\Hi}_1}$ may be absent ($= \{ 0 \}$).  \\   

\ni {\bf 6.  Implementation} \\

\ni {\it 1. Theory} \\
\i As remarked before, calculus is continuous; calculation is discrete.  Data is discrete.  Our given data on past curves (at times $1, \cdots, n$ say -- the distinction between being given the infinite past or a finite part of it is important theoretically but evaporates at the implementation stage) necessarily consists of finitely many points on each, interpolated or smoothed so as to give a curve (continuous, say), with whatever degree of smoothness the statistician chooses (as with density estimation; see e.g. Silverman [Sil]), by whatever means the statistician chooses -- splines, wavelets etc.  The prediction process then consists of the inevitable three steps: \\
(i) Discretisation of the data (now curves), into $d$-vectors, for some $d < \infty$.  For choice of $d$, see e.g. Li and Hsing [LiH]. \\ 
(ii) Prediction, using e.g. the multidimensional version of the Levinson-Durbin algorithm, as in [Bin2]; see below. \\
(iii) Interpolation or smoothing of this set of predicted values at time $n+1$ to give the predicted curve at time $n+1$. \\
See e.g. Aue et al. [AueNH] for a use of the Karhunen-Lo\`eve expansion here. \\
\i As one will see from the above, the practical problems involved at the implementation stage are largely {\it numerical}.  This is a familiar phenomenon; see e.g. [BinS2] for a different setting (random fields) where the writer (with Symons) recently encountered such things. \\
{\it The Levinson-Durbin algorithm} \\ 
\i For this classical algorithm for the computation of the best linear predictor based on the last $n$ data points, see e.g. Brockwell and Davis [BroD, \S 5.2], [Bin1] in the scalar case, [Bin2] in the matrix case.  It is worth noting that numerical improvements have been made here (the {\it split Levinson algorithm}); see Delsarte and Genin [DelG]. \\

\ni {\it 2. Numerics} \\
\i The three-step procedure above thus amounts substantially to\\
(i) {\it discretization} of the infinite-dimensional process (random curves) $Y$ in $(CR)$ to random $d$-vectors for suitably chosen $d$; \\
(ii){\it prediction} by finite-dimensional methods [Bin2] (e.g., the Levinson-Durbin algorithm); \\
(iii) {\it smoothing} (e.g., spline interpolation with a roughness penalty) to return to the infinite-dimensional setting. \\
See also Hyndman and Shang [HynS1], [HynS2]. \\
\i Use of the Karhunen-Lo`eve expansion here [MarR, \S 5.3], [GinN, Th. 2.6.10] (as in [AueNH] above) effectively specialises to the Gaussian case.  Such methods can be extremely effective (for background see e.g. [BinS1], [BinS2] and the references there); how appropriate they are in a particular case will depend on context. \\  
{\it Kernel methods} \\
\i Kernel methods, commonly used in machine learning, have recently been advocated for functional prediction by Hashimoto et al. [Has1], [Has2]. \\  

\ni {\bf 7. Complements} \\

\ni {\it 1. The deterministic case} \\
\i When the Szeg\H{o} condition $(Sz)$ fails (e.g., when the density is absent -- $F_a = 0$ in the notation of \S 3), the entire process is `bad', and consists entirely of `echoes of the remote past'.  In the scalar case (for simplicity): the variance ${\s}_n^2$ of the best linear predictor based on the last $n$ readings decreases to ${\s}^2 = 0$.  The interesting question of how fast was addressed long ago by  Rosenblatt [Ros], and more recently by Babayan, Ginovyan and Taqqu [BabGT].  \\

\ni {\it 2. Model spaces}  \\
\i Again in the scalar case first for simplicity: the (unilateral, forward) {\it shift} $S$ (time $n \mapsto n+1$) is represented on the Hardy space $H^2$ by 
$$
Sf = z f \qquad (f = f(z) \in H^2).
$$  
This has adjoint the backward shift $S^*$:
$$
S^*f = \frac{f - f(0)}{z}.
$$
This follows from the Taylor series for $f(z) = \sum_0^{\infty} a_n z^n$ in each case.  Note that $f \in H^2$ and $a = (a_n) \in {\ell}_2(\mathbb{N})$   are equivalent, and give the `Hardy-Hilbert space' (see de Branges [dBra], Martinez-Avendano and Rosenthal [MartR] for monograph treatments). \\ 
\i By the Beurling(-Lax-Halmos) theorem, the invariant subspaces of $S$ are $uH^2$ for the inner functions $u$.  Similarly, the invariant subspaces for $S^*$ are the orthogonal complements of these, written
$$
{\K}_u := (u H^2)^{\perp}.
$$
These are called {\it model spaces} (`model spaces are the invariant subspaces of the backward shift').  To explain the terminology, we quote [GarMR, 105] `The term {\it model space} originates in the theory of model operators, developed by Sz.-Nagy and Foia\c s, where it is shown that certain types of Hilbert-space contractions are unitarily equivalent to the compressions of the unilateral shift to a model space.  This underscores the importance of model spaces in developing concrete, function-theoretic realizations of abstract Hilbert space operators.'  For background, see the two classic sources [SzNFBK] (`first edition [SzNF]') and Nikolskii [Nik1] cited in \S 5, plus Rosenblum and Rovnyak [RosR] and the two-volume [Nik2], and the two recent treatments by Garcia, Mashreghi and Ross [GarMR] cited above (see also [GarR]) and Agler, McCarthy and Young [AglMY].  (Note that different authors use the term compression differently [GarMR, Remark 9.2].)  See \S 5 for references specific to the vector-valued case $H^2(\Hi)$, our main concern here, to which we add [Nik2, \S 4.8.8] and [NikV]. \\

\ni {\it 3. Compressions and dilations} \\
\i For $A, B$ operators on spaces ${\cal A} \subset {\cal B}$, $B$ is a {\it dilation} of $A$ if
$$
A^n = pr \ B^n \quad \hbox{for all} \ n \in \N,
$$
where $pr$ is projection [SzNBK, p.10].  Then $A$ is a {\it compression} of $B$ ([GarMR, Def. 9.1]; [SzNBK] does not use the term compression). \\
\i The basic result here is the {\it Sz.-Nagy-Foia\c s dilation theorem}: if  $H$ is a Hilbert space and $T$ a linear contraction on it, there exists a larger Hilbert space $\Hi$ and a unitary operator $U$ on $\Hi$ with $U$ a dilation of $T$ and 
$$
{{\Hi}} = {\bigvee} (U^n H: n \in \mathbb{N})
$$
(such a dilation is called {\it minimal}).  For the extensive theory here, see e.g. [SzNFBK, Ch. I], [Nik1, Introductory Lecture, Lecture III], [GarMR, Ch. 9].  \\
   
\ni {\it 4. Harmonizability}  \\
\i This concept, due to Lo\`eve [Lo\`e] in 1948, addresses the need to relax the strong assumption of stationarity, which one cannot expect to hold exactly in practice.   It has been studied and extended by Karhunen, Cram\'er and others.  As the covariance function now needs two arguments, $s$ and $t$ say, rather than one, $s-t$, the spectral measure now becomes a `bimeasure'.  For details and references, see e.g. [Rao1,2], [Kak1-4].  To summarise [Rao2, 292]:  \\
strongly harmonisable $\subset$ weakly harmonizable $\subset$ Karhunen $\subset$ Cram\'er. \\
\i Note [Kak3, Th. 2] that a process is weakly harmonizable iff it has a stationary dilation.  Thus weak harmonizability is the broadest context in which we can hope to bring the powerful tools available in the stationary case to bear.  Beyond that, one is more in the realm of the Kalman filter and its extensions, where the dynamics dominate, one predicts using only the immediate past and the present, and what matters is speed and accuracy of reaction (e.g., control of manned spacecraft, or mortar-locating radar). \\

\ni {\it 5. Banach spaces and beyond} \\
\i Such extensions were addressed briefly in [BinM, \S 5.5].  We refer there for further detail, and to e.g. Dette et al. [DetKA], Chobanyan and Weron [ChoW], Weron [Wer], Miamee and Salehi [MiaS] and Klotz and Riedel [KloR].  For a Banach-space version of the Wold decomposition, see [FauH]. \\
 
\ni {\it 6. Operator-valued processes}  \\
\i For an approach via dilation theory and operator models, see Makagon and Salehi [MakS, \S 2] and the references cited there.  The vector- and operator-valued cases are developed together in [RosR, Ch. 4-6].    \\ 
  
\ni {\it 7. The multivariate case} \\
\i See e.g. [Kak3, \S 6] for a multivariate vector-valued treatment. \\

\ni {\it 8. High- and infinite-dimensional probability and statistics}. \\
\i For high-dimensional treatments, see e.g. [Ver] for probability, [Wai] for statistics.  For infinite-dimensional statistics, see e.g. [Gin\'e and Nickl [GinN]. \\

\ni {\it 9. Specifically infinite-dimensional phenomena}. \\
\i We have been concerned here with vector-valued Hardy space theory (inner-outer function factorization, etc.), in finitely many dimensions [Bin2], or infinitely many, our main concern here.  There are results in the area which do {\it not} extend to the infinite-dimensional case.  For example, see Treil [Tre1], [Tre2] for the (operator) corona problem. \\

\ni {\it 10. Cram\'er-Karhunen-Lo\`eve expansion} \\
\i The Karhunen-Lo\`eve expansion ([MarR, \S 5.3], [GinN, Th. 2.6.10]) may be combined with the Cram\'er representation to give a dynamicized form of the expansion, which has been used for modelling of functional time series.  See e.g. Antoniadis and Sapatinas [AntS] (who considered El Ni\~ no), [AntPS], Panaretos and Tavakoli [PanT1], [PanT2].  \\
 
\ni {\it 11. Functional data analysis; prediction; change-points} \\
\i For some recent developments in these areas, see e.g. Aue, Norinho and H\"ormann [AueNH], Dette, Kokot and Aue [DetKA]. \\

\ni {\it 12.  Smoothness of functions} \\
\i Our data are functions, drawn from function spaces, often Hilbert spaces $\Hi$.  These may (and typically do) have some smoothness properties.  One of the most basic is {\it continuity of point evaluation}, $x \mapsto f(x)$ 
($f \in \Hi$).  By the Riesz representation theorem, this is the condition for $\Hi$ to be a {\it reproducing-kernel Hilbert space} (RKHS): to have a {\it reproducing kernel} $k(.,.)$ such that, with $k_x(.) := k(.,x)$,
$$
(f,k_x) = (f(.), k(.,x))_{\Hi} = f(x) \qquad (f \in {\Hi}).             
$$
Such spaces are common in the Hardy-space setting above.  For example, the Hardy-Hilbert space $H^2$ is a RKHS with kernel
$$
k_{\lambda}(z) = 1/(1 - \overline{\lambda} z) \qquad (\lambda, z \in \D)
$$
[GarMR, Prop. 3.3], and so are the model spaces ${\K}_u$ above, with kernel (with the $u$ here to be understood in the notation)
$$
k_{\lambda} (z) = \frac{1 - \overline{u(\lambda}) u(z)}{1 - \overline{\lambda} z} \qquad (\lambda, z \in \D)
$$
[GarMR, \S 5.5].  For the extensive theory of RKHS, and applications to probability and statistics, see e.g. Berlinet and Thomas-Agnan [BerTA].\\

\ni {\it 13.  Hankel operators and Nehari sequences}.  \\
\i Toeplitz operators occur frequently in the above; they have many links with Hankel operators (for which see e.g. Peller [Pel]).  The {\it Nehari problem} [Pel, Ch. 5] is: given a sequence $\g = ({\g}_n)$ (${\g}_n \in \C$), find a function $\phi$ in the unit ball of $L_{\infty}$ with
$$
{\g}_n = \int e^{in \t} \phi(\t) d\t/2 \pi \quad (n = 1,2, \cdots).
$$
{\it Nehari's theorem} is that such a solution exists iff the Hankel matrix  $({\g}_{m+n})$ of $\g$ acts as a contraction on ${\ell}_2$.  There is more than one solution (the {\it indeterminate case}) iff $\g = ({\g}_n)$ is the negatively indexed Fourier coefficients of the phase factor of some function in $H_2$; in this case $\g$ is called a {\it Nehari sequence} (compare the determinate and indeterminate cases in the moment problem; see e.g. [Bin3] and the references there). \\
\i  Nehari sequences also occur in prediction theory, for example in connection with the condition of {\it complete non-determinacy} for time series (see e.g. [KasB1]).  Study of problems of Nehari type has led (among other things) to extensions of the strong Szeg\H{o} and Baxter theorems in the scalar case [KasB1, \S 5] and to matrix forms of Baxter's theorem [KasB2,3].    \\ 
\\    

\ni {\bf 8.  Questions} \\

\i We close with some questions arising from the work surveyed above. \\
{\it Q1}.  Does the work of [DamPS] (\S 3) on Verblunsky's theorem extend to infinitely many dimensions? \\
{\it Q2}.  Does the work of [DerHKT] (\S 3) on Szeg\H{o}'s theorem extend to infinitely many dimensions? \\
{\it Q3}.  Is there an infinite-dimensional version of Szeg\H{o}'s theory of orthogonal polynomials on the unit circle (OPUC, [Sim1], [Bin1]; MOPUC, [Bin2])?  In this regard, see [GarMR, Example 8.3].  This involves the {\it Szeg\H{o}} (or {\it Toeplitz}) {\it conjugation}, which is familiar from the {\it Szeg\H{o} recursion} of OPUC and MOPUC. \\
{\it Q4}.  Can the assumption of essentially bounded spectral density in [Pow], \S 4, be relaxed (or dropped)? \\
{\it Q5}.  Can the assumption of having a scalar multiple in [SzNFBK, V.6,7], \S 4, be relaxed (or dropped)? \\
\\ 

\ni {\bf References} \\

\ni [AglMY] J. Agler, J. E. McCarthy and N. Young, {\sl Operator analysis: Hilbert space methods in complex analysis}.  Cambridge Tracts Math. {\bf 219}, Cambridge University Press, 2020. \\
\ni [AntPS] A. Antoniadis, E. Paroditis and T. Sapatinas, A functional wavelet-kernel approach for time series prediction.  {\sl J. Roy. Stat. Soc. B} {\bf 68} (2006), 837-857. \\
\ni [AntS]  A. Antoniadis and T. Sapatinas, Wavelet methods for continuous-time prediction using Hilbert-valued autoregressive processes.  {\sl J. Multivariate Analysis} {\bf 87} (2003), 133-158.  \\
\ni [AueNH] A. Aue, D. D. Norinho and S. H\"ormann, On the prediction of stationary functional time series.  {\sl J. Amer. Stat. Soc.} {\bf 110} (2015), 378-392. \\
\ni [AuevD] A. Aue and A. van Delft, Testing for stationarity of functional time series in the frequency domain.  {\sl Ann. Statist.} {\bf 48} (2020), 2505-2547; arXiv:1701.01741. \\
\ni [BabGT] N. K. Babayan, M. S. Ginovyan and M. S. Taqqu, Extensions of Rosenblatt's results on the asymptotic behaviour of the prediction error for deterministic stochastic sequences.  Special Issue: Murray Rpsenblatt Memorial,  {\sl J. Time Series Analysis} {\bf 42}(5-6) (2021), 622-652; arXiv:2006.00430. \\
\ni [BerTA] A. Berlinet and C. Thomas-Agnan, {\sl Reproducing kernel Hilbert spaces with applications in probability and statistics}.  With a preface by Persi Diaconis.  Kluwer, 2004.    \\
\ni [Beu] A. Beurling, On two problems concerning linear transformations in Hilbert space.  {\sl Acta Math.} {\bf 81} (1948), 239-255 (reprinted in {\sl Collected works of Arne Beurling}, Vol. {\bf 2}, {\sl Harmonic analysis} (ed. L. Carleson et al.), Birkh\"auser, 1989). \\ 
\ni [Bin1] N. H. Bingham, Szeg\H{o}'s theorem and its probabilistic descendants.  {\sl Probability Surveys} {\bf 9} (2012), 287-324.\\
\ni [Bin2] N. H. Bingham, Multivariate prediction and matrix Szeg\H{o} theory.  {\sl Probability Surveys} {\bf 9} (2012), 325-339. \\
\ni [Bin3] N. H. Bingham,  The life, work and legacy of P. L. Chebyshev.  {\it Proceedings: Bicentennial Conference on P. L. Chebyshev (1821-1894)} (ed. A. N. Shiryaev), {\sl Th. Probab. Appl.} {\bf 66}(4) (2021), to appear. \\
\ni [BinM] N. H. Bingham and Badr Missaoui, Aspects of prediction.  {\sl J. Applied Probability} {\bf 51A} (2014), 189-201. \\
\ni [BinS1] N. H. Bingham and Tasmin L. Symons, Gaussian random fields on sphere and sphere cross line.  {\sl Stochastic Proc. Appl.} (Larry Shepp Memorial Issue); arXiv:1812.02103; {\tt https://doi.org/10.1016/j.spa.2019.08.007}.
\ni [BinS2] N. H. Bingham and Tasmin L. Symons,  Aspects of random fields.  {\sl Theory of Probability and Mathematical Statistics} (Special Issue in honour of M. I. Yadrenko, ed. A. Olenko), to appear. \\
\ni [BleL] D. P. Blecher and L. E. Labuschagne, Applications of the Fuglede-Kadison determinant: Szeg\H{o}'s theorem and outers for noncommutative $H^p$.  {\sl Trans. Amer. Math. Soc.} {\bf 360} (2008), 6131-6147. \\
\ni [dBra] L. de Branges, {\sl Square-summable power series}.  Springer, Monographs in Math., 2010 (1st ed., with J. Rovnyak, Holt, Rinehart and Winston, 1966). \\
\ni [BroD] P. J. Brockwell and R. A. Davis, {\sl Time series: Theory and methods}, 2nd ed., Springer, 1991 (1st ed. 1987). \\ 
\ni [ChoW] S. A. Chobanyan and A. Weron, Banach space valued stationry processes and their linear prediction.  {\sl Dissertationes Math.} {\bf 125}, 1975. \\
\ni [Cra1] H. Cram\'er, On the theory of stationary random processes.  {\sl Ann. Math.} {\bf 41} (1940), 215-230 (reprinted in {\sl Collected Works of Harald Cram\'er, Volume II}, 925-940, Springer, 1994). \\
\ni [Cra2] H. Cram\'er, On harmonic analysis in certain function spaces.  {\sl Ark. Mat. Astr. Fys.} {\bf 28B} (1942), 1-7 ({\sl Works II}, 941-947). \\
\ni [Cra3] H. Cram\'er, A contribution to the theory of stochastic processes.  {\sl Proc. Second Berkeley Symposium Mat. Stat. Prob} (ed. J. Neyman) 329-339, U. California Press, 1951 ({\sl Works II}, 992-1002). \\
\ni [CraL] H. Cram\'er and R. Leadbetter, {\sl Stationary and related stochastic processes}.  Wiley, 1967. \\
\ni [DamPS] D. Damanik, A. Pushnitski and B. Simon, The analytic theory of matrix orthogonal polynomials. {\sl Surveys in Approximation Theory} {\bf 4} (2008), 1-85; arXiv:0711:2703. \\
\ni [vDelE] A. van Delft and M. Eichler, A note on Herglotz's theorem for time series on function spaces.  {\sl Stoch. Proc. Appl.} {\bf 130} (2020), 3687-3710.\\
\ni [DelG]  P. Delsarte and Y. Genin, The split Levinson algorithm.  {\sl IEEE Trans. Acoust. Speech Signal Proc} ASSP-{\bf 34} (1986), 470-478.  \\
\ni [DelGK] P. Delsarte, Y. Genin and Y. G. Kamp, Orthogonal polynomial matrices on the unit circle.  {\sl IEEE Trans. Ciruits and Systems} CAS-{\bf 25} (1978), 149-160. \\
\ni [DerHKT] M. Derevyagin, O. Holtz, S. Khrushchev and M. Tyaglov, Szeg\H{o}'s theorem for matrix orthogonal polynomials. {\sl J. Approx. Th.} {\bf 164} (2012), 1238-1261; arXiv:1104.4999. \\
\ni [DetKA] H. Dette, K. Kokot and A. Aue, Functional data analysis in the Banach space of continuous functions.  {\sl Ann. Stat.} {\bf 48} (2020), 1168-1192. \\
\ni [DunS] N. Dunford and J. T. Schwartz, {\sl Linear operators, II: Spectral theory, Self-adjoint operators in Hilbert space}, Wiley, 1963. \\
\ni [Dur] P. L. Duren, {\sl Theory of} $H^p$ {\sl spaces}.  Academic Press, 1970.\\
\ni [DymK] H. Dym and D. P. Kimsey, CMV matrices, a matrix version of Baxter's theorem, scattering and de Branges spaces.  {\sl EMS Surveys Math. Sci.} {\bf 3} (2016), 105p. \\
\ni [FauH] G. D. Faulkner and J. E. Honeycutt, Orthogonal decompositions of isometries in a Banach space.  {\sl Proc. Amer. Math. Soc.} {\bf 69} (1978), 125-128.   \\
\ni [GarMR] S. R. Garcia, J. Mashreghi and W. T. Ross, {\sl Introduction to model spaces and their operators}.  Cambridge Stud. Adv. Math. {\bf 148}, Cambridge University Press, 2016. \\
\ni [GarR] S. R. Garcia and W. T. Ross, Model spaces: A survey.  {\sl Invariant subspaces of the shift operator}, 197-245.  {\sl Contemp. Math.} {\bf 638}, Amer. Math. Soc., 2015. \\
\ni [Garn] J. B. Garnett, {\sl Bounded analytic functions}.  Academic Press, 1981 (revised 1st ed., Grad. Texts Math. {\bf 236}, Springer, 2007). \\
\ni [GinN] E. Gin\'e and R. Nickl, {\sl Mathematical foundations of infinite-dimensional statistical models}.  Cambridge University Press, 2016. \\ 
\ni [GreSz] U. Grenander and G. Szeg\H{o}, {\sl Toeplitz forms and their applications}.  U. California Press, 1958. \\
\ni [Gret] N. E. Gretsky, Operator valued Gramians and inner products in vector valued function spaces.  {\sl Ann. Mat. Pura Appl.} {\bf 101} (1976), 337-354. \\
\ni [Hal] P. R. Halmos, Shifts on Hilbert spaces.  {\sl J. Reine Angew. Math.} {\bf 208} (1961), 102-112. \\
\ni [Has1]  Y. Hashimoto, I. Ishikawa, M. Ikeda, F. Kumura and Y. Kawahara, Kernel mean embeddings of von Neumann-algebra-valued measures.  arXiv:2007.14698. \\ 
\ni [Has2] Y. Hashimoto, I. Ishikawa, M. Ikeda, F. Kumura, T. Katsura and Y. Kawahara, Reproducing kernel Hilbert $C^{\ast}$-module and kernel mean embedding.  arXiv:2101.11410.  \\
\ni [Hof] K. Hoffman, {\sl Banach spaces of analytic functions}.  Prentice-Hall, 1962. \\ 
\ni [HorJ] R. A. Horn and C. R. Johnson, {\sl Matix analysis}, Cambridge University Press, 1985. \\
\ni [HorK] L. H\'orvath and P. Kokoszka, {\sl Inference for functional data with applications}.  Springer, 2012. \\
\ni [Horo] J. Horowitz, Gaussian random measures.  {\sl Stoch. Proc. Appl.} {\bf 22} (1986), 129-133. \\
\ni [HynS1] R. J. Hyndman and H. L. Shang, Forecasting functional time series (with discussion).  {\sl J. Korean Math. Soc.} {\bf 38} (2009), 199-221. \\
\ni [HynS2] R. J. Hyndman and H. L. Shang, The {\tt ftsa} package for R. \\ {\tt https.//cran.r-project.org/web/packages/ftsa/ftsa.pdf}. \\
\ni [Kak1] Y. Kakihara, A classification of vector harmonisable processes.  {\sl Stoch. Anal. Appl.} {\bf 10} (1992), 277-311. \\
\ni [Kak2] Y. Kakihara, Vector harmonizable processes: Wold and Cram\'er decompositions.  {\sl Stoch. Anal. Appl.} {\bf 13} (1995), 531-541. \\
\ni [Kak3] Y. Kakihara, The Kolmogorov isomorphism theorem and extensions to some nonstationary processes.  {\sl Stochastic processes: Theory and methods, Handbook of Statistics } {\bf 19} (ed. D. N. Shanbhag and C. R. Rao), 443-470, North-Holland, 2001. \\
\ni [Kak4] Y. Kakihara, Spectral domains of vector harmonizable processes.  {\sl J. Statistical Planning and Inference} {\bf 100} (2002), 93-108. \\
\ni [KalM] G. Kallianpur and V. Mandrekar, Spectral theory of stationary $H$-valued processes.  {\sl J. Multiv. Analysis} {\bf 1} (1971), 1-16. \\
\ni [KasB1] Y. Kasahara and N. H. Bingham, Verblunsky coefficients and Nehari sequences.  {\sl Trans. Amer. Math. Soc.} {\bf 366} (2014), 1363-1378.  \\
\ni [KasB2] Y. Kasahara and N. H. Bingham, Coefficient stripping in the matrix Nehari problem.  {\sl J. Approx. Th.} {\bf 220} (2017), 1-11.   \\
\ni [KasB3] Y. Kasahara and N. H. Bingham, Matricial Baxter's theorem with  Nehari sequence.  {\sl Math. Nachrichten} {\bf 291} (2018), 2590-2598.  \\
\ni [KloR] L. Klotz and M. Riedel, Spectral representation and extrapolation of stationary random processes on linear spaces.  {\sl Prob. Math. Stat.} {\bf 21} (2001), 179-197. \\
\ni [Lan] E. C. Lance, {\sl Hilbert} $C^*$-{\sl modules: A toolkit for operator algebraists}.  London Math. Soc. Lect. Notes Ser. {\bf 210}, Cambridge Univ. Press, 1995.    \\
\ni [Lax] P. D. Lax, Translation invariant spaces.  {\sl Acta Math.} {\bf 101} (1959), 163-178.\\
\ni [LiH] Y. Li and T. Hsing, Deciding the dimension of effective dimension reduction space for functional and high-dimensional data.  {\sl Ann. Stat.} {\bf 38} (2010), 3028-3062. \\
\ni [Lo\`e] M. Lo\`eve, Fonctions al\'eatoires du second ordre.  Appendix to P. L\'evy, {\sl Processus stochastiques et mouvement brownien}, Gauthier-Villars, 1948. \\
\ni [Mac] G. W. Mackey, {\sl Unitary group representations in physics, probability and number theory}.  Benjamin, 1978. \\
\ni [MakS] A. Makagon and H. Salehi, Notes on infinite dimensional stationary sequences.  Prob. Th. on Vector Spaces IV, {\sl Lecture Notes in Math.} {\bf 1391} (1989), 200-238. \\
\ni [ManS] V. Mandrekar and H. Salehi, The square-integrability of operator-valued functions with respect to a non-negative operator-valued measure and Kolmogorov's isomorphism theorem.  {\sl Indiana U. Math. J.} {\bf 20} (1972), 545-563. \\
\ni [MarR] Marcus, M. B. and Rosen, J., {\sl Markov processes, Gaussian processes and local times}.  Cambridge University Press, 2006. \\
\ni [MartR] R. A. Martinez-Avendano and P. Rosenthal, {\sl An introduction to operators on the Hardy-Hilbert space}.  Grad. Texts Math. {\bf 237}, Springer, 2007. \\
\ni [Mas1] P. Masani, Orthogonally scattered measures.  {\sl Adv. Math.} {\bf 2} (1968), 61-117. \\
\ni [Mas2] P. Masani, Comments on the prediction-theoretic papers.  {\sl Norbert Wiener: Coll. Works Vol. III} (ed. P. Masani), MIT Press, 1981, p.276-306. \\
\ni [MiaS] A. G. Miamee and H. Salehi, The factorization problem for positive operator-valued functions in a Banach space and regularity of Banach space valued stationary stochastic processes.  {\sl Sanhy\=a A} {\bf 39} (1977), 211-222. \\ 
\ni [Nik1] N. K. Nikolskii, {\sl Treatise on the shift operator: Spectral function theory}. Grundl. math. Wiss. {\bf 273}, Springer, 1986.\\
\ni [Nik2] N. K. Nikolskii, {\sl Operators, functions and systems: an easy reading.  Volume 1: Hardy, Hankel and Toeplitz; Volume 2: Model operators and systems}.  Math. Surveys and Monographs {\bf 92, 93}, Amer. Math. Soc., 2002. \\
\ni [NikV] N. K. Nikolskii and V. I. Vasyunin, Notes on two function models.  {\sl The Bieberbach conjecture: Proceedings of the Symposium on the Occasion of the Proof}, 113-141.  Math. Surv. Monog. {\bf 21}, Amer. Math. Soc., 1986. \\
\ni [PanT1] V. Panaretos and S. Tavakoli, Fourier analysis of stationary time series in function space.  {\sl Ann. Stat.} {\bf 41} (2013), 568-603.\\
\ni [PanT2] V. Panaretos and S. Tavakoli, Cram\'er-Karhunen-Lo\`eve representation and harmonic principal components analysis of functional time series.  {\sl Stoch. Proc. Appl.} {\bf 123} (2013), 2779-2807. \\
\ni [Pay] R. Payen, Fonctions al\'eatoires du second ordre \`a valeurs dans un espace de Hilbert.  {\sl Ann. Inst. H. Poincar\'e Prob. Stat.} {\bf 3} (1967), 323-396. \\
\ni [Pel] V. V. Peller, {\sl Hankel operators and their applications}.  Springer, 2003. \\
\ni [Pon] L. S. Pontryagin, {\sl Topological groups}, 2nd ed., Gordon and Breach, 1966. \\ 
\ni [Pow] S. C. Power, Spectral characterization of the Wold-Zasuhin decomposition and prediction-error operator.  {\sl Math. Proc. Cambridge Phil. Soc.} {\bf 110} (1991), 559-567. \\
\ni [RamS1] J. O. Ramsay and B. W. Silverman, {\sl Functional data analysis}, 2nd ed., Springer, 2005 (1st ed. 1997). \\
\ni [RamS2] J. O. Ramsay and B. W. Silverman, {\sl Applied functional data analysis: Methods and case studies}, Springer, 2002. \\
\ni [Rao1] M. M. Rao, Harmonizable processes: structure theory.  {\sl Enseignement Math.} {\bf 28} (1982), 295-351. \\
\ni [Rao2] M. M. Rao, Harmonizable, Cram\'er, and Karhunen classes of process.  {\sl Handbook of statistics} {\bf 5}: {\sl Time series in the time domain} (ed. E. J. Hannan et al.), North-Holland, 1985, 279-310. \\
\ni [RieN] F. Riesz and B. Sz.-Nagy, {\sl Le\c cons d'analyse fonctionelle}, 2nd ed., Akad\'emiai Kiad\'o, 1953 (1st ed. 1952). \\
\ni [Ros] M. Rosenblatt, Some purely deterministic processes.  {\sl J. Math. Mech.} {\bf 6} (1957), 801-810 (reprinted in {\sl Selected works of Murray Rosenblatt} (ed. R. A. Davis et al.), Springer, 2011, 124-133). \\
\ni [RosR] M. Rosenblum and J. Rovnyak, {\sl Hardy classes and operator theory}.  Oxford University Press, 1985; Dover, 1997. \\
\ni [Rud1] W. Rudin, {\sl Real and complex analysis}.  McGraw-Hill, 3rd ed., 1987 (1st ed. 1966, 2nd ed. 1974). \\
\ni [Rud2] W. Rudin, {\sl Functional analysis}, 2nd ed., McGraw-Hill, 1991 (1st ed. 1973). \\
\ni [Sil] B. W. Silverman, {\sl Density estimation for statistics and data analysis}.  Monogr. Stat. Appl. Prob. {\bf 26}, Chapman \& Hall, 1986. \\
\ni [Sim1] B. Simon, {\sl Orthogonal polynomials on the unit circle.  Part 1: Classical theory}.
AMS Colloquium Publications {\bf 54.1}, American Math. Soc., Providence RI, 2005. \\
\ni [Sim2] B. Simon, {\sl Orthogonal polynomials on the unit circle.  Part 2: Spectral theory}.  AMS Colloquium Publications {\bf 54.2}, American Math. Soc., Providence RI, 2005. \\
\ni [Sim3] B. Simon, {\sl Szeg\H{o}'s theorem and its descendants: Spectral theory for} $L^2$ {\sl perturbations of orthogonal polynomials}.  Princeton University Press, 2011.\\
\ni [Sto] M. H. Stone, {\sl Linear transformations in Hilbert space and their applications to analysis}.  Amer. Math. Soc. Colloq. Publ. {\bf XV}, AMS, 1932. \\
\ni [SzNFBK] B. Sz.-Nagy, C. Foia\c s, H. Bercovici and L. K\'erchy, {\sl Harmonic analysis of operators on Hilbert space}, 2nd ed., Springer Universitext, 2010 (1st ed., Sz.-Nagy and Foia\c s, North-Holland, 1970).\\
\ni [Tre1] S. R. Treil, Angles between co-invariant subspaces and the operator corona problem.  The Sz\"okefalvi-Nagy problem.  {\sl Dokl. Akad. Nauk SSSR} {\bf 302} (1988), 1063-1068. \\
\ni [Tre2] S. R. Treil, An operator corona theorem.  {\sl Indiana U. Math. J.} {\bf 53} (2004), 1763-1780. \\
\ni [Ver] R. Vershynin, {\sl High-dimensional probability.  An introduction with applications in data science}.  Cambridge University Press, 2018. \\
\ni [Wai] M. J. Wainwright, {\sl High-dimensional statistics: A non-asymptotic viewpoint}.  Cambridge University Press, 2019. \\
\ni [Wer] A. Weron, Prediction theory in Banach spaces.  Proc. Winter Sch. Prob. Karpacz.  {\sl Lecture Notes in Math.} {\bf 472} (1975), 207-228.\\
\ni [Whi] P. Whittle, The analysis of multiple time series.  {\sl J. Roy. Stat. Soc. B} {\bf 15} (1953), 125-139. \\

\ni Mathematics Department, Imperial College, London SW7 2EZ, UK; n.bingham@ic.ac.uk \\ 

\end{document}